\documentclass[10pt]{article}
\usepackage{graphicx} 

\usepackage{amsmath, amsfonts, amsthm, amssymb}

\usepackage{bm}
\usepackage{todonotes}
\usepackage[paper=a4paper,total={17cm, 27cm}]{geometry}

\usepackage{subcaption}

\usepackage{xcolor}

\newcommand{\pf}{\varphi}
\newcommand{\Div}{\nabla\cdot}
\newcommand{\strain}[1][u]{{\bm \varepsilon}(\bm #1)}
\newcommand{\out}{\mathrm{out}}
\newcommand{\inn}{\mathrm{in}}
\newcommand{\bn}{\mathbf{n}}
\newcommand{\lnabla}{\nabla_{\Gamma}}
\newcommand{\Sym}{\mathcal S}
\newcommand{\Cf}{\mathbb C}

\theoremstyle{remark}

\usepackage{authblk}

\title{Sharp-Interface Limit of the Cahn-Hilliard-Biot Equations}

\author[1]{Erlend Storvik\footnote{Corresponding author: erlend.storvik@hvl.no}}
\author[1]{Carina Bringedal}

\affil[1]{Western Norway University of Applied Sciences, Department of Computer science, Electrical engineering and Mathematical sciences}

\date{}

\begin{document}

\maketitle

\subsection*{Abstract}
In this letter, we derive the sharp-interface limit of the Cahn-Hilliard-Biot equations using formal matched asymptotic expansions. We find that in each sub-domain, the quasi-static Biot equations are obtained with domain-specific material parameters. Moreover, across the interface, material displacement and pore pressure are continuous, while volumetric fluid content and normal stress are balanced. By utilizing the energy of the system, the phase-field potential is shown to be influenced by the curvature, along with contributions from both flow and elasticity at the interface. The normal velocity of the interface is proportional to the jump in normal derivative of the phase-field potential across the interface. Finally, we present a numerical experiment that demonstrates how the location of each phase evolves consistently as the diffuse-interface width parameter becomes smaller; only the width of the diffuse interface changes.

\section{Introduction and model equations}
The Cahn-Hilliard-Biot system is a diffuse-interface model coupling the Cahn-Hilliard equations for phase evolution to the quasi-static Biot equations for poroelasticity. A thorough derivation of the model can be found in \cite{storvik2022}. Recently, the well-posedness of the model was studied in \cite{agh2025, riethmuller2024, fritz2023well}, and solutions strategies have been discussed in \cite{bf2025, storvik2024}. However, the sharp-interface limit of the model equations is yet to be investigated in literature.

Let $\Omega \subset \mathbb{R}^d$ be a domain with a Lipschitz continuous boundary, $d \in \{2,3\}$ be the spatial dimension and $T_\mathrm{f} > 0$ be the final time. We consider the problem: Find $(\pf, \mu, \bm u, p, \theta)$ for $({ x},t)\in \Omega \times [0,T_\mathrm{f}]$, with $\varphi$ being the phase-field variable, $\mu$ the chemical potential, $\bm u$ the displacement, $p$ the pore pressure, and $\theta$ the volumetric fluid content, such that

\begin{subequations}
\begin{align}
\partial_t \varphi - \Div (m \nabla \mu) &= R, \label{eq:ch1}\\
\mu + {\gamma}\left(\ell\Delta\varphi-\frac{1}{\ell}\Psi'(\varphi)\right) - \delta_\varphi\mathcal{E}_\mathrm{e}(\varphi, \bm u) - \delta_\varphi \mathcal{E}_\mathrm{f}(\varphi, \bm u, \theta) &= 0, \label{eq:ch2}\\
 -\nabla\cdot\big(\mathbb{C}(\varphi)\left(\bm\varepsilon\left(\bm u\right)-\mathcal{T}(\varphi)\right)-\alpha(\varphi)p\bm I\big) &= {\bm f}, \label{eq:elasticity}\\ 
\partial_t\theta - \Div ({\kappa(\pf)\nabla p}) &= S_\mathrm{f}, \label{eq:flow}\\
 p - M(\varphi)\left(\theta - \alpha(\varphi)\nabla\cdot \bm u\right) &=0,\label{eq:ptheta}
\end{align}
\label{eq:model}
\end{subequations}
where 
\begin{equation}\label{eq:diff_elasticity}\delta_\varphi\mathcal{E}_\mathrm{e}(\varphi, \bm u) = -\mathcal{T}'\!(\varphi)\!:\!\mathbb{C}(\varphi)\big(\bm\varepsilon\left(\bm u\right)-\mathcal{T}(\varphi)\big)+\frac{1}{2} \big(\strain - \mathcal{T}(\varphi)\big)\!:\!\mathbb{C}'\!(\varphi)\big(\strain - \mathcal{T}(\varphi)\big),
\end{equation}
and 
\begin{equation}\label{eq:diff_fluid}\delta_\varphi\mathcal{E}_{\mathrm{f}}(\varphi,\bm u, \theta) = \frac{M'(\varphi)}{2}\left(\theta - \alpha(\pf)\nabla\cdot \bm u\right)^2-M(\pf)\left(\theta-\alpha(\pf)\nabla\cdot \bm u\right)\alpha'(\varphi)\nabla\cdot \bm u,
\end{equation}
accompanied with appropriate initial and boundary conditions. 
Here, $m$ is the chemical mobility, $R$ is a reaction term,  $\gamma > 0$ denotes the surface tension, $\ell > 0$ is a small parameter associated with the width of the regularization layer,  $\Psi(\varphi) := \varphi^2\left(1-\varphi\right)^2$ is a double-well potential penalizing deviations from the pure phases, having local minima at $\pf = 1$ and $\pf = 0$,  $\mathbb{C}(\varphi)$ is the stiffness tensor, $\strain := \frac{1}{2}\left(\nabla \bm u + \nabla\bm u^\top\right)$ is the linearized strain tensor, $\mathcal{T}(\pf):= \xi(\pf - \bar{\pf})\bm I$ accounts for swelling effects, $\alpha(\pf)$ is the Biot-Willis coupling coefficient, $\bm f$ accounts for body forces, $M(\pf)$ is the compressibility coefficient, $S_f$ is a source term, and $\kappa(\varphi)$ is the permeability. The terms in \eqref{eq:diff_elasticity} and \eqref{eq:diff_fluid} correspond to the variational derivatives of the potential elastic energy and the hydraulic energy of the system with respect to the phase-field.

The goal of this letter is to establish the sharp-interface limit of the Cahn-Hilliard-Biot equations \eqref{eq:model}. By employing matched asymptotic expansions, we will in the following find that in each of the phases, we recover the quasi-static Biot equations with desired material parameters, while across the interface the material displacement and pore pressure are continuous, and stress and volumetric fluid content are balanced. Furthermore, we will exploit the free energy of the system and derive a condition for the chemical potential at the interface. The derivation of the sharp-interface limit follows a similar procedure to the ones used for the Cahn-Hilliard-Navier-Stokes equations in \cite{agg2012} and for phase-field modeling of precipitation and dissolution in \cite{bringedal2020}.

\section{The sharp-interface limit}
We use the method of formal matched asymptotic expansions to investigate the behavior of the Cahn-Hilliard-Biot equations away from and close to the interface.

\subsection{Inner and outer expansions}

The outer expansions consider the solution in regions away from the interface. We mark all variables in this view with a superscript $y^{\out}$, and have the following expansions $\sum_{i=0}^\infty \ell^i y^\out_i(x,t)$, 
 for $y\in\{\pf,\mu,\bm u,p,\theta\}$.

The inner expansions consider the solution in an interfacial region where the transition between the two phases takes place. We make a change in variables to a local coordinate system near the interface, where $s$ parametrizes the location on the interface (where $\pf=0.5$), and $r$ is the signed distance away from the interface. We then use the scaled distance $z=\frac{r}{\ell}$. Therefore, the inner expansions are $\sum_{i=0}^\infty \ell^iy_i^\inn(z,s,t)$
 for $y\in\{\pf,\mu,\bm u,p,\theta\}$.
 Derivatives need to be rewritten according to the coordinate transformation. We refer to e.g.~\cite{agg2012,caginalp1988dynamics} for a justification of the following rules:
\begin{align}
\partial_t f &= -\frac{1}{\ell}v_n \partial_z f^\inn + \text{h.o.t.}, \ \  \nabla_x f = \frac{1}{\ell}\partial_z f^\inn \bm n+\lnabla f^\inn+ \text{h.o.t.}, \ \ \nabla^2_x f = \frac{1}{\ell^2}\partial_z^2f^\inn -\frac{H}{\ell}\partial_zf^\inn  +  \text{h.o.t.},\\
\nabla_x\cdot \bm f &= \frac{1}{\ell}\partial_z \bm f^\inn\cdot \bm n+\lnabla\cdot  \bm f^\inn+ \text{h.o.t.}, \quad 
\nabla_x \bm f = \frac{1}{\ell}\partial_z\bm f^\inn \otimes \bn + \lnabla\bm f^\inn+ \text{h.o.t.},\\
\bm\varepsilon(\bm f) &= \frac{1}{2}\Big[\frac{1}{\ell}\partial_z\bm f^\inn\otimes \bn + \frac{1}{\ell}\bn \otimes \partial_z \bm f^\inn + \lnabla \bm f^\inn + (\lnabla\bm f^\inn)^T\Big]+ \text{h.o.t.}, \label{eq:innerstrain}
\end{align}
where $v_n$ is the velocity of the interface, $H$ the mean curvature, $\lnabla$ the surface gradient, $\bn$ the normal vector of the interface, and h.o.t. stands for higher order terms. 
Moreover, we have the matching conditions
\begin{align}
    \lim_{z\to\pm \infty} f_0^\inn(z, s) &= f_0^\out(x\pm),\label{eq:match0}\\
    \lim_{z\to\pm\infty} \partial_z f_1^\inn(z, s) &= \nabla f_0^\out(x\pm)\cdot \bn.\label{eq:match1}
\end{align}

In the following, we first consider the leading order of the outer expansions in Section \ref{sec:outer}, which will give us model equations within the two sub-domains. Afterwards we consider the leading order of the inner expansions in Section \ref{sec:inner}, which gives some properties of the variables at the interface, and afterwards the next order terms of the inner expansions in Section \ref{sec:inner2} give coupling conditions between the two domains.

 \subsection{Leading order of the outer expansions }\label{sec:outer}

We now insert the outer expansions into the Cahn-Hilliard-Biot equations 
\eqref{eq:ch1}--\eqref{eq:ptheta}. Considering \eqref{eq:ch2}, we observe that the leading order term, of order $\ell^{-1}$, will be $-\gamma\Psi'(\varphi^\out_0)$, which admit solutions $\pf_0 = 0, 0.5$ or $1$, where $\pf_0=0,1$ are stable.
We now denote by $\Omega^0$ and $\Omega^1$ the regions in which $\pf^\out_0 = 0$ and $\pf^\out_0 = 1$, respectively. By doing so, the leading order terms of equations \eqref{eq:elasticity}--\eqref{eq:ptheta} yield the standard quasi-static Biot equations in each subdomain, with their respective material parameters:
\begin{align}
 -\nabla\cdot\big(\mathbb{C}_{i}\left(\bm\varepsilon\left(\bm u_0^\out\right)-\mathcal{T}_i\right)-\alpha_ip_0^\out\bm I\big) &= {\bm f} \quad x\in \Omega^i,\label{eq:elasticity_outer}\\ 
\partial_t\theta_0^\out - \nabla\cdot(\kappa_i\nabla p_0^\out) &= S_\mathrm{f} \quad x\in \Omega^i, \label{eq:flow_outer}\\
p_0^\out - M_i(\theta_0^\out-\alpha_i\nabla\cdot \bm u_0^\out) &= 0 \quad x\in \Omega^i,\label{eq:darcyflow_outer}
\end{align}
where $\mathbb{C}_i=\mathbb{C}(i)$, $\mathcal{T}_i = \mathcal{T}(i)$, $\alpha_i=\alpha(i)$, $M_i=M(i)$ and $\kappa_i = \kappa(i)$ for $i=0,1$.

 \subsection{Leading order of the inner expansions }\label{sec:inner}
We now insert the inner expansions into equation \eqref{eq:ch1}--\eqref{eq:ptheta} and consider here the leading order terms. 
As we insert the inner expansions into \eqref{eq:flow}, 
the dominating term of order $\frac{1}{\ell^2}$ is $\partial_z(\kappa(\pf^\inn_0)\partial_z p^\inn_0 \bn)\cdot \bn =0$. 
Because $\bn$ is independent on $z$ we have 
$\partial_z(\kappa(\pf^\inn_0)\partial_z p^\inn_0)=0$
which by anti differentiation and matching condition \eqref{eq:match0} corresponds to $\partial_z p_0^\inn = 0.$ 
Integrating from $-\infty$ to $\infty$ and matching condition \eqref{eq:match0} gives that
\begin{equation}\label{eq:p0hopp}
    \left[p_0^\out\right]_-^+=0.
\end{equation}

The dominating term of equation \eqref{eq:ptheta} of order $\frac{1}{\ell}$ is $-M(\pf_0^\inn)\alpha(\pf_0^\inn)\partial_z \bm u_0^\inn\cdot \bn = 0\Leftrightarrow \partial_z \bm u_0^\inn\cdot \bn = 0$,
which by integration from $-\infty$ to $\infty$ and matching condition \eqref{eq:match0} gives
\begin{equation}
    \left[\bm u_0^\out\cdot \bn\right]_-^+=0
\end{equation}

For \eqref{eq:elasticity}, we take some extra care with rewriting the derivatives, as this will be relevant for the next order term. The calculation is analogous to \cite{agg2012} and the linearized strain tensor is rewritten corresponding to \eqref{eq:innerstrain}.
We introduce $\Sym(A) = \frac{1}{2}(A+A^T)$ for a square matrix $A$. Then
\begin{align}
    \nabla\cdot(\mathbb{C}(\pf)(\strain-\mathcal{T}(\pf)) = &\frac{1}{\ell^2}\partial_z(\Cf(\pf^\inn)\Sym(\partial_z\bm u^\inn\otimes\bn)\bn) + \lnabla\cdot(\Cf(\pf^\inn)(\Sym(\lnabla \bm u^\inn)-\mathcal{T}(\pf^\inn)))\nonumber\\ &+ \frac{1}{\ell} \Big[\partial_z(\Cf(\pf^\inn)(\Sym(\lnabla \bm u^\inn)-\mathcal{T}(\pf^\inn))\bn)  + \lnabla\cdot(\Cf(\pf^\inn)\Sym(\partial_z \bm u^\inn\otimes \bn) \Big].\label{eq:derivative}
\end{align}
Therefore, as we insert the inner expansions in \eqref{eq:elasticity}, the leading order term of order $\frac{1}{\ell^2}$ is $\partial_z(\Cf(\pf_0^\inn)\Sym(\partial_z\bm u^\inn_0\otimes\bn)\bn)=0$. 
Writing out the symmetric part, and using properties of outer product along with $\partial_z \bm u_0^\inn\cdot \bn = 0$, we obtain $\partial_z(\Cf(\pf^\inn_0)\partial_z\bm u_0^\inn )=0$. 
Integrating from $-\infty$ to $\infty$ gives us that $\bm u_0^\inn$ has to be constant in $z$, and using matching condition \eqref{eq:match0} this means
\begin{equation}\label{eq:elastinnerlimit}
    [\bm u_0^\out]^+_-=0,
\end{equation}
namely that the entire displacement is continuous across the interface.

For \eqref{eq:ch2}, we take some extra care on how to rewrite the derivatives in the term coming from \eqref{eq:diff_elasticity}. We have that
\begin{align*}
    \delta_\varphi\mathcal{E}_\mathrm{e}(\varphi, \bm u) = &-\frac{1}{\ell}\xi\bm I\!:\!\Cf(\pf)(\Sym(\partial_z \bm u^\inn\otimes \bn) - \xi\bm I\!:\!\Cf(\pf)(\Sym(\lnabla\bm u^\inn)-\mathcal T(\pf^\inn)) \\ &+\frac{1}{\ell^2}\frac{1}{2}\Sym(\partial_z \bm u^\inn\otimes \bn)\!:\!\Cf'(\pf^\inn)\Sym(\partial_z \bm u^\inn\otimes \bn) +
    \frac{1}{\ell}(\Sym(\lnabla\bm u^\inn) -\mathcal T(\pf^\inn))\!:\!\Cf'(\pf^\inn)\Sym(\partial_z\bm u^\inn\otimes\bn) \\
      &+ \frac{1}{2}(\Sym(\lnabla\bm u^\inn) -\mathcal T(\pf^\inn))\!:\!\Cf'(\pf^\inn)(\Sym(\lnabla\bm u^\inn) -\mathcal T(\pf^\inn)).
\end{align*}
The dominating terms when inserting the inner expansions into \eqref{eq:ch2} is of order $\frac{1}{\ell}$ and gives that
\begin{equation}
    \partial_z^2\pf^\inn_0 - \Psi'(\pf^\inn_0) = 0. \label{eq:ch2leading}
\end{equation}
Note that there is a term of order $\frac{1}{\ell^2}$ coming from $\delta_\varphi\mathcal{E}_\mathrm{e}(\varphi, \bm u)$, but as it involves $\partial_z\bm u_0^\inn$, this term vanishes. Terms of order $\frac{1}{\ell}$ coming from $\delta_\varphi\mathcal{E}_\mathrm{e}(\varphi, \bm u)$ and $\delta_\varphi\mathcal{E}_\mathrm{f}(\varphi, \bm u,\theta)$ are also zero as they depend on $\partial_z\bm u_0^\inn$.
From \eqref{eq:match0} we have that $\lim_{z\to\infty}\pf^\inn_0 = 1$ and $\lim_{z\to-\infty}\pf^\inn_0 = 0$. 
Further, we let $\pf^\inn_0(z=0)=0.5$ as this should correspond to the location of the evolving interface. With our choice of double-well potential, \eqref{eq:ch2leading} has the solution $\pf^\inn_0(z) = \frac{1}{1+e^{-\sqrt{2}z}} = \frac{1}{2}(1+\tanh({\frac{z}{\sqrt{2}}}))$,
which gives us the form of the diffuse transition, and more importantly, that $\pf^\inn_0$ only depends on $z$.

The leading order term of \eqref{eq:ch1} is of order $\frac{1}{\ell^2}$ and is $\partial_z(m\partial_z\mu^\inn_{0}\bn) \cdot \bn = 0$, 
which means that $\mu_{0}^\inn$ is constant in $z$. 
By \eqref{eq:match0}, we obtain
\begin{equation}\label{eq:ch1innerlimit}
    [\mu_{0}^\out]^+_- = 0.
\end{equation}

\subsection{Next order of the inner expansions}\label{sec:inner2}
Here, we inspect equation \eqref{eq:ch1}--\eqref{eq:ptheta} in the inner expansions of the next order of $\ell$.

\subsubsection{Interfacial phase-field flux balance}
The next order terms of \eqref{eq:ch1} are of order $\frac{1}{\ell}$ and are
$-v_n\partial_z\pf_0^\inn - \partial_z(m\partial_z\mu_{1}^\inn) =0$ 
since $\mu_0^\inn$ is constant in $z$. 
By integrating and using matching condition \eqref{eq:match1}, we obtain
\begin{equation}\label{eq:fluxbalance}
    v_n = -m[\nabla\mu_{0}^\out\cdot\bn]^+_- .
\end{equation}

\subsubsection{Interfacial flow balance}
The next order terms of \eqref{eq:flow} are of order $\frac{1}{\ell}$ and are $-v_n\partial_z\theta_0^\inn -\partial_z(\kappa(\pf_0^\inn)(\partial_zp_1^\inn\bn +\lnabla p_0^\inn))\cdot\bn = 0$. 
 By integrating in $z$ from $-\infty$ to $+\infty$ and using matching conditions \eqref{eq:match0} and \eqref{eq:match1}, we obtain
\begin{equation}\label{eq:flowbalance}
    v_n[\theta_0^\out ]^+_- + [\kappa_i\nabla p_0^\out\cdot \bn]^+_-=0,
\end{equation}
where $\kappa_i$, $i=0,1$, is the permeability at the respective side of the interface. 
The next order terms of \eqref{eq:ptheta} are
\begin{equation}\label{eq:p_equality}
p_0^\inn - M(\pf_0^\inn)(\theta_0^\inn -\alpha(\pf_0^\inn)(\partial_z\bm u_1^\inn\cdot\bn + \lnabla\cdot\bm u_0^\inn)) = 0.
\end{equation}
Letting $z\to\pm\infty$ gives us that
\begin{equation*}
p_0^\out|_i - M_i(\theta_0^\out|_i - \alpha_i\nabla\cdot \bm u_0^\out|_i) =0,
\end{equation*}
where $\alpha_i$ and $M_i$, $i=0,1$, are the material properties found at the respective side of the interface., and $|_i$ indicates that the variable should be evaluated at the interface but at the respective side of the interface.

\subsubsection{Interfacial stress balance}
For the next order terms of \eqref{eq:elasticity} we make use of the rewritten derivatives from \eqref{eq:derivative}. We therefore have that
\begin{align}\label{eq:stressdiff}
    -\partial_z\left(\Cf(\pf_0^\inn)\left(\Sym(\partial_z\bm u_1^\inn\otimes\bn)- \Sym(\lnabla\bm u_0^\inn)-\mathcal{T}(\pf_0^\inn)\right)-\alpha(\pf_0^\inn)p_0^\inn\bm I\right)\bn = 0.
\end{align}
Terms that involve $\partial_z\bm u_0^\inn$ have already been removed as they are zero. The above terms are integrated in $z$ from $-\infty$ to $+\infty$ and by \eqref{eq:match0} and \eqref{eq:match1} we obtain
\begin{equation}\label{eq:elasticbalance}
    [(\Cf_i(\bm \varepsilon(\bm u_0^\out) -\mathcal{T}_i)-\alpha_ip_0^\out\bm I)\bn]^+_- =0,
\end{equation}
that is, the jump in the normal stress is zero. Here, $\mathcal T_i$ accounts for the swelling effects, with either $\pf_0^\out=0,1$ inserted according to the side of the interface.

\subsubsection{A generalized Gibbs-Thomson equation}
Next order terms of \eqref{eq:ch2}, that are of order $\ell^0$, are
\begin{align}
    \nonumber&\mu_{0}^\inn +\gamma \partial_z^2\pf_1^\inn - \gamma H\partial_z\pf_0^\inn -\gamma\Psi''(\pf_0^\inn)\pf_1^\inn \\ &-\frac{1}{2}\big(\Sym(\lnabla \bm u_0^\inn)+\Sym(\partial_z \bm u_1^\inn\otimes \bn) - \mathcal{T}(\pf_0^\inn)\big)\!:\!\Cf'(\pf_0^\inn)\big(\Sym(\lnabla \bm u_0^\inn)+\Sym(\partial_z \bm u_1^\inn\otimes \bn) - \mathcal{T}(\pf_0^\inn)\big) \nonumber\\
    &\nonumber + \xi\bm I\!:\!\Cf(\pf_0^\inn)(\Sym(\partial_z \bm u_1^\inn\otimes\bn)+\Sym(\lnabla \bm u_0^\inn)-\mathcal T(\pf_0^\inn))- \frac{M'(\pf_0^\inn )}{2}\left(\theta_0^\inn-\alpha(\pf_0^\inn)(\partial_z\bm u_1^\inn \cdot \bn + \lnabla \cdot\bm u_0^\inn)\right)^2\\
    & + M(\pf_0^\inn)\left(\theta_0^\inn-\alpha(\pf_0^\inn)(\partial_z\bm u_1^\inn \cdot \bn + \lnabla \cdot\bm u_0^\inn)\right) \alpha'(\pf_0^\inn)(\partial_z\bm u_1^\inn \cdot \bn + \lnabla \cdot\bm u_0^\inn)  = 0. \label{eq:next_order_mu}
\end{align}
We now introduce the notation $\bm e\left(\varphi_0^\inn, \bm u_0^\inn, \bm u_1^\inn\right):= \Sym(\lnabla \bm u_0^\inn)+\Sym(\partial_z \bm u_1^\inn\otimes \bn) - \mathcal{T}(\pf_0^\inn)$
and the function
\begin{align}
K\left(\varphi_0^\inn, \bm u_0^\inn, \bm u_1^\inn, \theta_0^\inn\right):=&\frac{1}{2}\bm e\left(\varphi_0^\inn, \bm u_0^\inn, \bm u_1^\inn\right)\!:\!\Cf(\pf_0^\inn)\bm e\left(\varphi_0^\inn, \bm u_0^\inn, \bm u_1^\inn\right)\nonumber\\&+\frac{M(\pf_0^\inn)}{2} \left(\theta_0^\inn - \alpha(\pf_0^\inn)\left(\partial_z\bm u_1^\inn \cdot \bn + \lnabla \cdot\bm u_0^\inn\right)\right)^2.\label{eq:gibbs-thomson-energy}
\end{align}
Differentiating $K$ with respect to $z$, and recalling that $\bm u_0^\inn$ is independent of $z$, gives
\begin{equation*}
    \partial_z K\left(\varphi_0^\inn, \bm u_0^\inn, \bm u_1^\inn, \theta_0^\inn\right)=K_{\delta \pf_0^\inn}\left(\varphi_0^\inn, \bm u_0^\inn, \bm u_1^\inn, \theta_0^\inn\right) + K_{\delta \bm u_1^\inn}\left(\varphi_0^\inn, \bm u_0^\inn, \bm u_1^\inn, \theta_0^\inn\right) + K_{\delta \theta_0^\inn}\left(\varphi_0^\inn, \bm u_0^\inn, \bm u_1^\inn, \theta_0^\inn\right),
\end{equation*}
where $K_{\delta \pf_0^\inn}$, $K_{\delta \bm u_1^\inn}$, $K_{\delta \theta_0^\inn}$ correspond to the derivatives of $K$ related to $\pf_0^\inn$, $\bm u_1^\inn$ and $\theta_0^\inn$, respectively.

We multiply \eqref{eq:next_order_mu} with $\partial_z\pf_0^\inn$ and integrate from $-\infty$ to $\infty$. Upon using the matching conditions and \eqref{eq:ch2leading} we obtain
\begin{align*}
    \mu_{0}^\out &= \frac{\gamma H}{3\sqrt{2}} + \int_{-\infty}^\infty K_{\delta \pf_0^\inn}\left(\varphi_0^\inn, \bm u_0^\inn, \bm u_1^\inn, \theta_0^\inn\right) dz.
\end{align*}
By investigating the integrals of $K_{\delta \bm u_1^\inn}$ and $K_{\delta \theta_0^\inn}$ with respect to $z$ from $-\infty$ to $\infty$, and using equation \eqref{eq:p_equality}, that $p_0^\inn$ does not depend on $z$, and the matching conditions, we obtain
\begin{align*}
    \int_{-\infty}^\infty K_{\delta \bm u_1^\inn}\left(\varphi_0^\inn, \bm u_0^\inn, \bm u_1^\inn, \theta_0^\inn\right) dz = \int_{-\infty}^\infty(\Sym(\partial_z^2\bm u_1\otimes \bn)\!:\!\Cf(\pf_0^\inn)\bm e\left(\varphi_o^\inn, \bm u_0^\inn, \bm u_1^\inn\right)
    - \alpha(\pf_0^\inn)p_0^\inn\bm I \bn \cdot\partial_z^2\bm u_1^\inn) dz,
\end{align*}
\begin{equation*}
\int_{-\infty}^\infty K_{\delta \theta_0^\inn}\left(\varphi_0^\inn, \bm u_0^\inn, \bm u_1^\inn, \theta_0^\inn\right) dz = \int_{-\infty}^\infty p_0^\inn\partial_z\theta_0^\inn dz = p_0^\inn\left[\theta_0^\out\right]_-^+=\left[p_0^\out\theta_0^\out\right]_-^+.
\end{equation*}
   Using properties of the symmetric operator $\Sym$ and the outer product together with integration by parts, \eqref{eq:stressdiff}, and matching conditions, give
   \begin{align*}
   \int_{-\infty}^\infty &K_{\delta \bm u_1^\inn}\left(\varphi_0^\inn, \bm u_0^\inn, \bm u_1^\inn, \theta_0^\inn\right) dz = -\int_{-\infty}^\infty\partial_z\bm u_1\cdot\partial_z\left(\Cf(\pf_0^\inn)\bm e\left(\varphi_o^\inn, \bm u_0^\inn, \bm u_1^\inn\right)-\alpha(\pf_0^\inn)p_0^\inn\bm I\right)\bn dz\\
    &+\lim_{k\to \infty} \left[\partial_z\bm u_1\cdot \left(\Cf(\pf_0^\inn)\bm e\left(\varphi_o^\inn, \bm u_0^\inn, \bm u_1^\inn\right)-\alpha(\pf_0^\inn)p_0^\inn\bm I\right)\bn\right]_{-k}^k = \left[\nabla \bm u_0^\out\bn\cdot\left(\mathbb{C}_i\left(\bm \varepsilon(\bm u_0^\out)-\mathcal{T}_i\right)-\alpha_ip_0^\out\bm I\right)\bn\right]_-^+.
   \end{align*}
Hence, by adding and subtracting appropriate terms, we obtain
\begin{align}\label{eq:ch2gt}
    \mu_{0}^\out &= \frac{\gamma H}{3\sqrt{2}} +\int_{-\infty}^\infty \partial_z K dz - \left[\nabla \bm u_0^\out\bn\cdot\left(\mathbb{C}_i\left(\bm \varepsilon(\bm u_0^\out)-\mathcal{T}_i\right)-\alpha_ip_0^\out\bm I\right)\bn\right]_-^+
    \ - \left[p_0^\out\theta_0^\out\right]_-^+ ,
\end{align}
where
\begin{equation*}
\int_{-\infty}^\infty \partial_z K dz = \left[\frac{1}{2}\left(\bm \varepsilon(\bm u_0^\out)-\mathcal{T}_i\right)\!:\!\mathbb{C}_i\left(\bm \varepsilon(\bm u_0^\out)-\mathcal{T}_i\right) + \frac{M_i}{2}\left(\theta_0^\out-\alpha_i\nabla\cdot\bm u_0^\out\right)^2 \right]_-^+.
\end{equation*}

\section{Numerical experiment}\label{sec:sim}
Here, we present a numerical experiment to investigate the behavior of the model \eqref{eq:model} as $\ell$ becomes smaller. For this we solve the equations \eqref{eq:ch1}--\eqref{eq:ptheta}, in pressure-flux formulation (see e.g., \cite{storvik2024}) with three different values of $\ell$, namely $\ell = 0.1, 0.05, 0.025$. The phase-field dependent material parameters are treated in the same way as in \cite{storvik2024} with the interpolation function $\pi(\varphi)$, and we use material parameters from Table~\ref{tab:param}. Lam\'e parameters ($G,\lambda$) for isotropic stiffness are applied. We impose the initial condition $\varphi = 0$, for $x<=0.5$ and $\varphi = 1$ for $x>0.5$ and zero conditions for $\bm u$ and $p$. Homogeneous Dirichlet boundary conditions are applied to the elasticity equation and Dirichlet pressure boundary conditions with $p=1$ are applied at the left boundary and $p=0$ at the right boundary. The equations are discretized in time and space similarly as in \cite{storvik2024} with a semi-implicit time-discretization and standard finite element spaces. Newton's method is used to solve the nonlinear discrete equations with absolute incremental tolerance of $1e-6$. The implementation was done using Fenicsx \cite{dolfinx} and the source code can be found at https://github.com/EStorvik/chb.

\begin{table}[h]
    \centering
    \begin{tabular}{c||c|c|c|c|c|c|c|c|c|c|c|c|c|c}
       Symbol &$\gamma$  &  $m$  & $\xi$ & $\tilde{\varphi}$&  $G_0$, $G_1$ & $\lambda_0$, $\lambda_1$ &$M_0$, $M_1$  &$\kappa_0$, $\kappa_1$ & $\alpha_0$, $\alpha_1$ &$dt$&$h$ & $R$ & $\bm f$ & $S_f$ \\
       Value& $1$ & $1$ & $0.1$ & $0.5$ & $100$, $1$& $20$, $0.1$ &$1$, $1$ &  $1$, $0.01$ &  $1$, $1$ &  $1e-3$ &  $\frac{\sqrt\ell}{3.2}$& $0$& $0$& $0$
    \end{tabular}
    \caption{Simulation parameters used in numerical experiment.}
    \label{tab:param}
\end{table}

\begin{figure}[h]
\centering
    \begin{subfigure}{0.075\textwidth}
        \includegraphics[trim = 0 -55 0 0, clip, width = \textwidth]{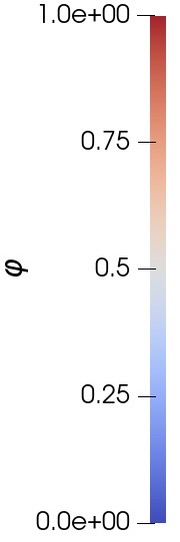}
    \end{subfigure}
    \begin{subfigure}{0.2\textwidth}
        \includegraphics[trim= 600 100 600 100,clip,width = \textwidth]{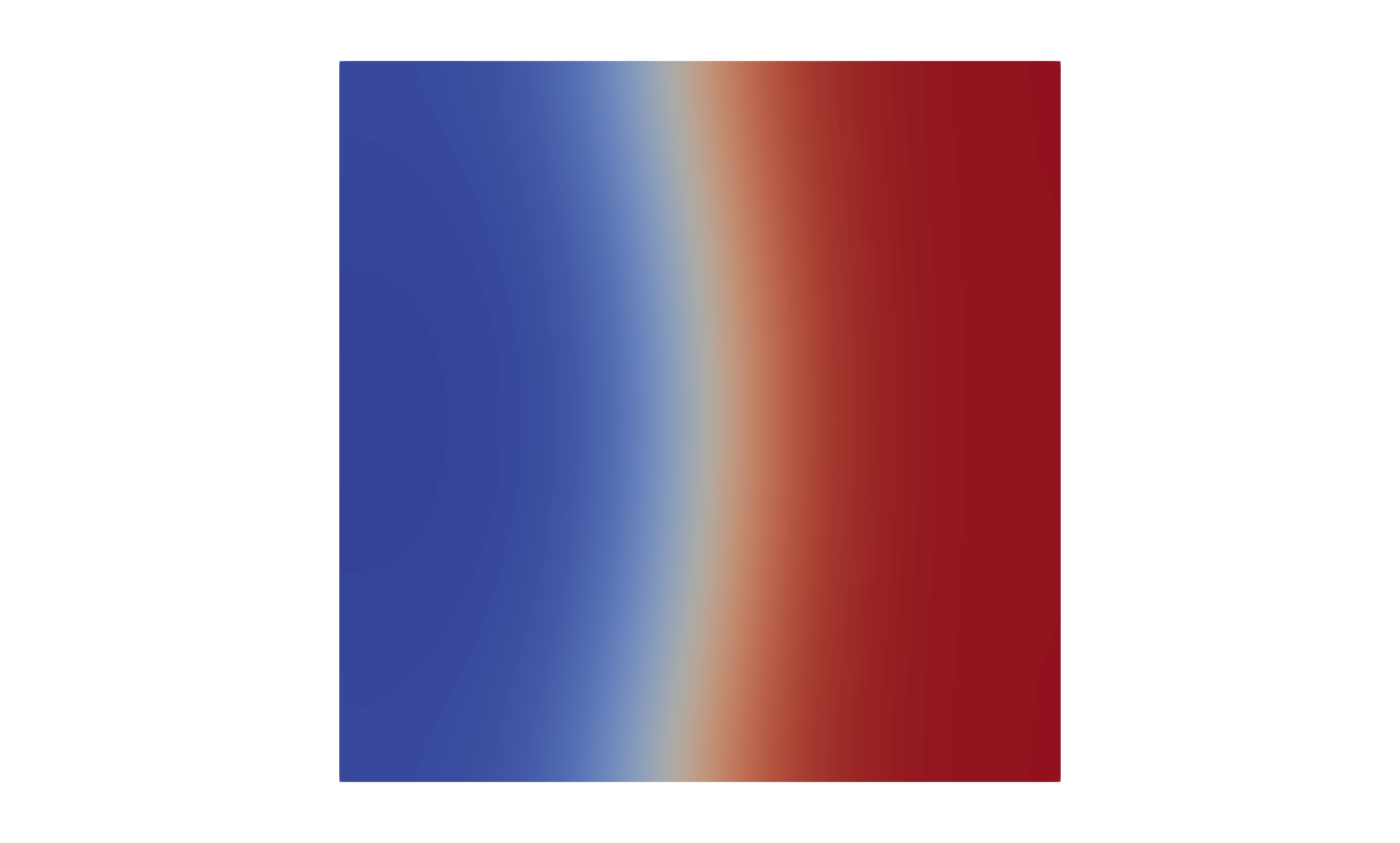}
        \caption{$\ell = 0.1$}
    \end{subfigure}
    \begin{subfigure}{0.2\textwidth}
        \includegraphics[trim= 600 100 600 100,clip,width = \textwidth]{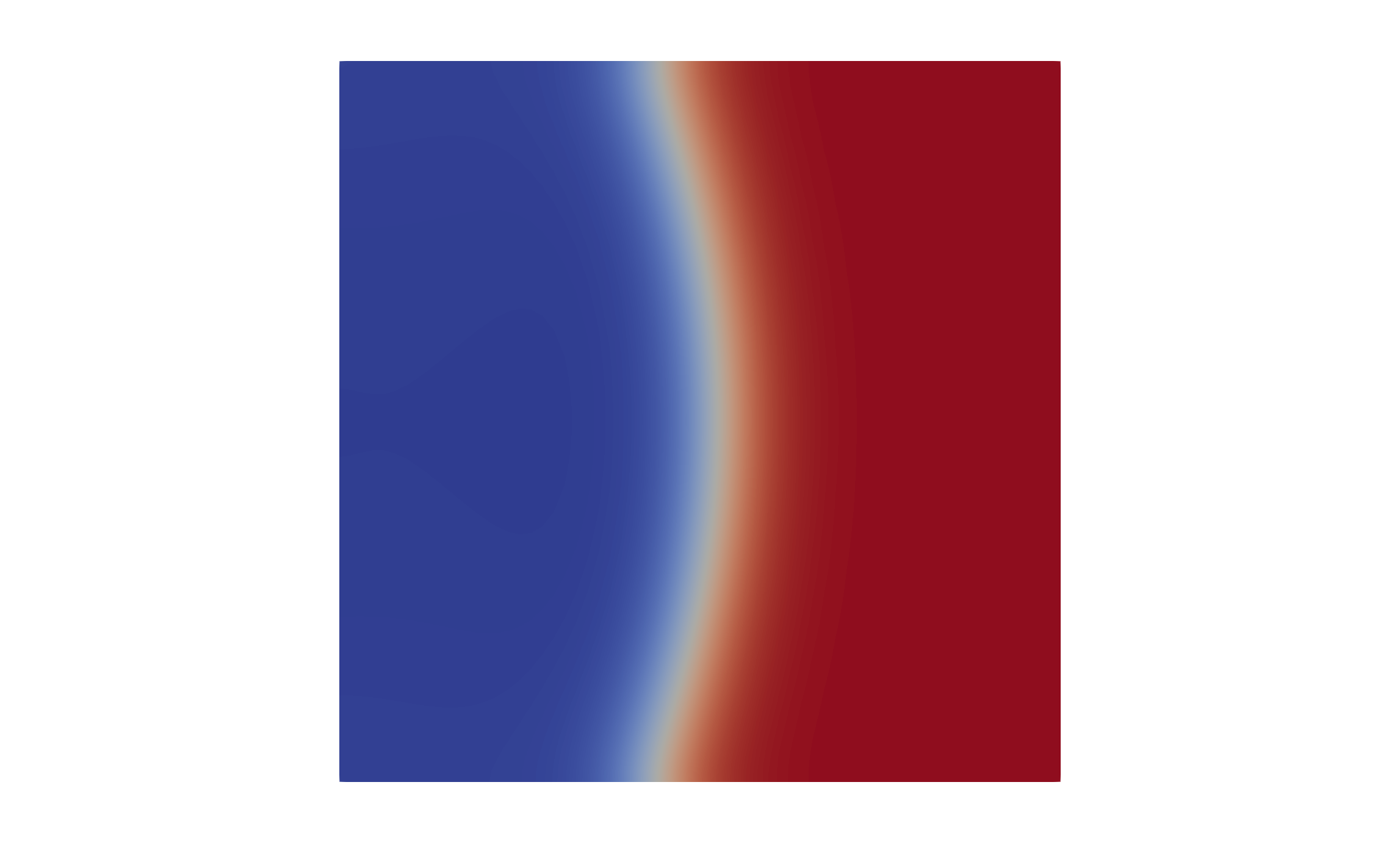}
        \caption{$\ell = 0.05$}
    \end{subfigure}
    \begin{subfigure}{0.2\textwidth}
        \includegraphics[trim= 600 100 600 100,clip,width = \textwidth]{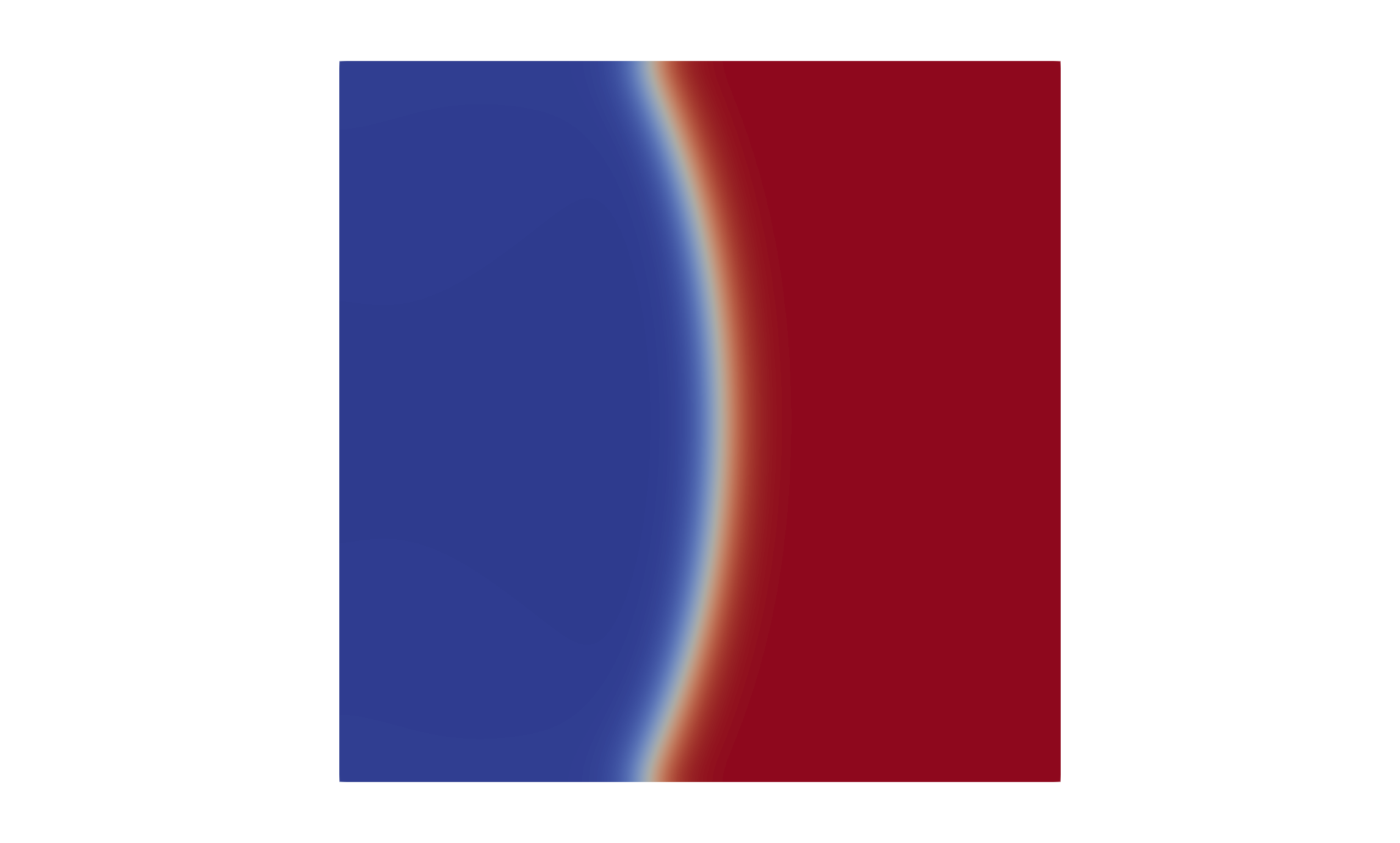}
        \caption{$\ell = 0.025$}
    \end{subfigure}
    \begin{subfigure}{0.295\textwidth}
        \includegraphics[width = \textwidth]{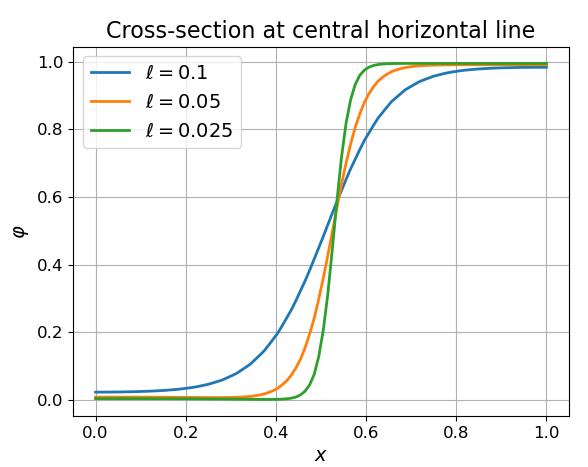}
        \caption{Cross-section of $\varphi$}
        \label{fig:cs}
    \end{subfigure}
       \caption{Comparison of the phase field at the end-time $t = 0.1$ for three different values of $\ell$, $0.1$, $0.05$ and $0.025$. Figure~\ref{fig:cs} shows the cross-section of the same phase fields ($\varphi$ at $t=0.1$) along the horizontal line going through the middle of the domain.}
       \label{fig:CHB}
\end{figure}

Simulation results are seen in Figure \ref{fig:CHB}. We observe from the experiment that the location of each phase evolves in a consistent manner for all three values of $\ell$ and that the interface sharpens as $\ell\to 0$, as is to be expected from the analysis.

\section{Conclusion}
In this letter, we investigated the sharp-interface limit of the Cahn-Hilliard-Biot equations and found that it represents quasi-static Biot equations \eqref{eq:elasticity_outer}-\eqref{eq:darcyflow_outer} in each subdomain, with coupling conditions \eqref{eq:fluxbalance}, \eqref{eq:flowbalance}, \eqref{eq:elasticbalance} and \eqref{eq:ch2gt}, and additionally \eqref{eq:p0hopp} and \eqref{eq:elastinnerlimit}. That is,
\begin{subequations}
\begin{align}
 -\nabla\cdot\big(\mathbb{C}_{i}\left(\bm\varepsilon\left(\bm u\right)-\mathcal{T}_i\right)-\alpha_ip\bm I\big) &= {\bm f} \quad x\in \Omega^i,\\ 
\partial_t\theta - \nabla\cdot(\kappa_i\nabla p) &= S_\mathrm{f} \quad x\in \Omega^i, \\
p- M_i(\theta-\alpha_i\nabla\cdot \bm u) &= 0 \quad x\in \Omega^i,
\end{align}
\end{subequations}
and with coupling conditions
\begin{subequations}
    \begin{align}
         v_n = -m[\nabla \mu\cdot\bn]^+_- &\quad x\in \Gamma,\\
         \mu = \frac{\gamma H}{3\sqrt{2}} + \left[\frac{1}{2}\left(\bm \varepsilon(\bm u)-\mathcal{T}_i\right)\!:\!\mathbb{C}_i\left(\bm \varepsilon(\bm u)-\mathcal{T}_i\right) + \frac{M_i}{2}\left(\theta-\alpha_i\nabla\cdot\bm u\right)^2 \right]^+_- \nonumber \\- \left[\nabla \bm u\bn\cdot\left(\mathbb{C}_i\left(\bm \varepsilon(\bm u)-\mathcal{T}_i\right)-\alpha_ip\bm I\right)\bn\right]_-^+
    \ - \left[p\theta\right]_-^+ &\quad x\in \Gamma,\\
         [(-\Cf_i(\bm \varepsilon(\bm u) -\mathcal{T}_i)+\alpha_ip)\bn]^+_- =0 &\quad x\in \Gamma,\\
         v_n[\theta ]^+_- + [\kappa_i\nabla p\cdot \bn]^+_-=0 &\quad x\in \Gamma,\\
         [\bm u]^+_- = 0 \text{ and } [p]^+_- = 0 &\quad x\in \Gamma,
    \end{align}
\end{subequations}
where $p|_i = M_i(\theta|_i - \alpha_i\nabla\cdot \bm u|_i) $.

Moreover, we provided a numerical experiment that validates that the location of each phase evolves in a consistent manner as $\ell$ becomes smaller, but the diffuse-interface width sharpens.
\bibliographystyle{unsrt}
\bibliography{bibliography}

\end{document}